\documentclass[11pt]{amsart}

\usepackage{amsmath, amssymb, amsthm}
\usepackage{enumitem}
\usepackage{graphicx}
\usepackage{hyperref}
\usepackage{mathtools}
\usepackage{subcaption}

\usepackage[utf8]{inputenc}
\usepackage[T1]{fontenc}

\title[A Family of Non-Periodic Tilings with Structural Regularity]{A Family of Non-Periodic Tilings, Describable Using Elementary Tools and Exhibiting a New Kind of Structural Regularity}
\author{Miki Imura}
\date{June 9, 2025}

\theoremstyle{definition}
\newtheorem{definition}{Definition}[section]

\theoremstyle{plain}
\newtheorem{theorem}[definition]{Theorem}
\newtheorem{proposition}[definition]{Proposition}
\newtheorem{lemma}[definition]{Lemma}
\newtheorem{observation}[definition]{Observation}

\graphicspath{{figures/}}
\setlist[itemize]{leftmargin=1.5em}

\begin{document}

\begin{abstract}
We present a construction of a family of non-periodic tilings using elementary tools such as modular arithmetic and vector geometry.
These tilings exhibit a distinct type of structural regularity, which we term \emph{modulo-staggered rotational symmetry}.
The construction is self-contained and does not rely on previous tiling theories or systems.
\end{abstract}

\maketitle

\section{Introduction}

\begin{figure}[htbp]
  \centering
  \begin{subfigure}[b]{0.45\textwidth}
    \centering
    \includegraphics[width=\textwidth]{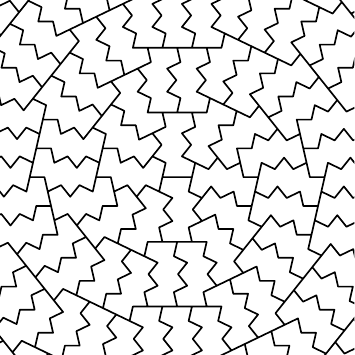}
    \caption{Without offset}
    \label{fig:no_offset}
  \end{subfigure}
  \hspace{0.03\textwidth}
  \begin{subfigure}[b]{0.45\textwidth}
    \centering
    \includegraphics[width=\textwidth]{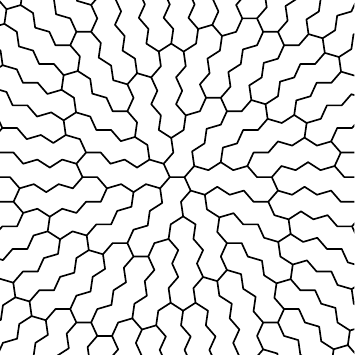}
    \caption{With offset}
    \label{fig:with_offset}
  \end{subfigure}
  \caption{Two variants of the Modulo Krinkle tiling.}
  \label{fig:intro}
\end{figure}

The study of tilings has led to a wide variety of patterns and mathematical techniques~\cite{grunbaum2016}.
Examples range from the classic Penrose tilings~\cite{penrose1979} to recent advances such as aperiodic monotiles~\cite{smith2024,smith2024chiral}.
In this paper, we propose a family of non-periodic tilings that can be described using only elementary tools.
While some previous studies~\cite{voderberg1936,gailiunas2000} also exhibit spiral-like structures, our construction relies on polygonal boundaries derived from modular arithmetic, yielding distinct patterns with a new kind of structural regularity (Figure~\ref{fig:intro}).

The structure of the paper is as follows.
In Section~\ref{sec:prelim}, we introduce basic notations and a framework for representing polygonal paths using a sequence and certain unit vectors.
In Section~\ref{sec:prototile}, we construct the prototile.
Section~\ref{sec:wedge} explains how to assemble wedge-shaped regions using the prototile.
In Section~\ref{sec:tiling}, we describe how to tile the plane using these regions.
Section~\ref{sec:proof} provides a formal proof of the validity of the construction.
Finally, Section~\ref{sec:symmetry} discusses its distinct structural regularity.

\section{Preliminaries}
\label{sec:prelim}

In this section, we introduce basic notations and a framework for constructing polygonal paths using a direction sequence and unit vectors placed at equal angular intervals around the circle.

We use $a \bmod b$ to denote the remainder of $a$ divided by $b$.
We denote sequences by $(a_j)$, indexed from $j = 0$ unless otherwise stated.

Let $m$ and $k$ be fixed positive integers such that $m < k$ and $\gcd(m, k) = 1$.
Let $t \ge 2$ be an integer.
We define $n$ according to the presence of offset:
\[
n = \begin{cases}
tk & \text{(without offset)}, \\
2(tk - m) & \text{(with offset)}.
\end{cases}
\]
By construction, we always have $2k \le n$, which will be used in subsequent arguments.

We define $\vec{v}_i = \left( \cos\frac{2\pi i}{n},\, \sin\frac{2\pi i}{n} \right)$ for $i = 0, 1, \dots, n-1$, representing unit vectors pointing in \emph{direction~$i$}.
Rotating $\vec{v}_i$ counterclockwise by $\frac{2\pi d}{n}$ yields $\vec{v}_{i + d}$.

Given a starting point $\vec{p} \in \mathbb{R}^2$ and a sequence of integers $(a_j)$ with $0 \le a_j < n$, we define a sequence of points $(\vec{p}_j)$ in $\mathbb{R}^2$ by
\[
\vec{p}_0 \coloneqq \vec{p}, \quad
\vec{p}_{j+1} \coloneqq \vec{p}_j + \vec{v}_{a_j} \text{ for } j \ge 0.
\]
We refer to such a sequence $(a_j)$ as a \emph{direction sequence}, and to the polygonal path formed by connecting the points $(\vec{p}_j)$ as the \emph{curve defined by $\vec{p}$ and $(a_j)$}.

\section{Construction of the Prototile}
\label{sec:prototile}

We now describe how to construct the prototile in our tiling system.
Let $(s_j)$ be a modular progression, defined by $s_j = jm \bmod k$, which consists of the remainders when successive multiples of $m$ are divided by $k$.
It is well known that the first $k$ terms of this sequence form a permutation of $0, 1, \dots, k-1$.

Define $(\ell_j)$ by taking these $k$ terms and appending $k$ to the end.
Then, define $(u_j)$ by swapping the first and last elements of $(\ell_j)$.

We define two polygonal paths from a common starting point~$P$:
the \emph{lower path} is the curve defined by $P$ and $(\ell_j)$, and
the \emph{upper path} is the curve defined by $P$ and $(u_j)$.

Observe that the two paths start at the same point~$P$, and since their direction sequences are permutations of each other, they end at the same location.
In addition, the directions appearing in the interior segments of each path are all strictly between $0$ and $k$.
Since $k \le \frac{n}{2}$ by construction, these directions span an angle strictly less than $180^\circ$\!.
Taking the inner product with the vector $\left( \cos\frac{k\pi}{n},\, \sin\frac{k\pi}{n} \right)$ shows that each step has a positive inner product with this vector, indicating a consistent orientation.
Since the interior segments of both paths run parallel and never intersect, the combined figure formed by the lower path followed by the reverse of the upper path defines a simple polygon, which serves as the prototile used in the tiling construction.

\begin{figure}[htbp]
  \centering
  \begin{subfigure}[b]{0.45\textwidth}
    \centering
    \includegraphics[width=\textwidth]{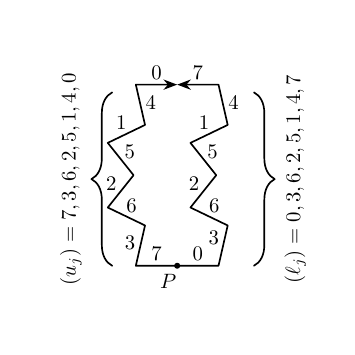}
    \caption{$n = 14$}
    \label{fig:prototile_no_offset}
  \end{subfigure}
  \hspace{0.03\textwidth}
  \begin{subfigure}[b]{0.45\textwidth}
    \centering
    \includegraphics[width=\textwidth]{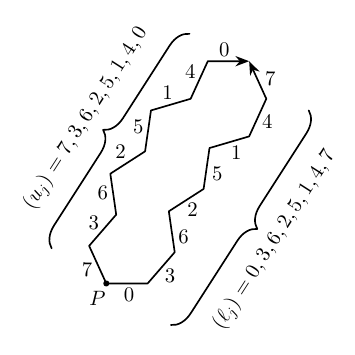}
    \caption{$n = 22$}
    \label{fig:prototile_with_offset}
  \end{subfigure}
  \caption{Two examples of the prototile constructed for $(m, k) = (3, 7)$, with different values of~$n$. Each edge is annotated with a direction label for clarity. When $n = 2k$, the interior angles at the starting point and ending point are exactly~$180^\circ$\!, causing the first and last edge pairs to appear as single straight segments of doubled length. However, we still treat them as two separate edges in our construction.}
  \label{fig:prototile_examples}
\end{figure}

Examples of this construction for $(m, k) = (3, 7)$, corresponding to the parameters used in the initial figures, are shown in Figure~\ref{fig:prototile_examples}.

To facilitate later discussion, we introduce terminology for parts of the prototile.

\begin{figure}[htbp]
  \centering
  \includegraphics[width=0.45\textwidth]{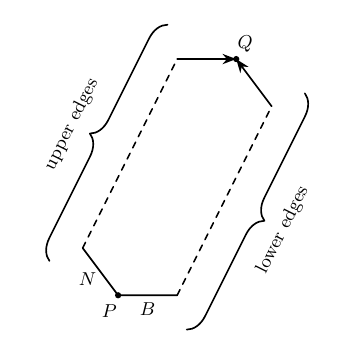}
  \caption{A conceptual diagram of the prototile, with the starting point~$P$, ending point~$Q$, base edge~$B$, and neighbor edge~$N$ labeled for reference.}
  \label{fig:prototile_concept}
\end{figure}

We call the common starting point of the lower and upper paths the \emph{starting point} of the tile, and their shared endpoint as the \emph{ending point}.
The first segment of the lower path is called the \emph{base edge}, and the first segment of the upper path is called the \emph{neighbor edge}, as it is adjacent to the base edge.
We refer to all segments of the lower path as the \emph{lower edges}, and those of the upper path as the \emph{upper edges}.
The last segment of the upper path, which is parallel to the base edge, is referred to as the \emph{opposite edge of the base edge}.
Similarly, the last segment of the lower path, which is parallel to the neighbor edge, is referred to as the \emph{opposite edge of the neighbor edge} (Figure~\ref{fig:prototile_concept}).

\section{Building a Wedge-Shaped Region}
\label{sec:wedge}

We next describe how to construct a wedge-shaped region---a key component in the full tiling---by placing congruent copies of the prototile in the same orientation.
While we do not provide a fully rigorous proof in the formal sense---such as defining a function that maps every point $(x, y)$ in the wedge to the unique tile it lies in---the construction clearly yields a regular tiling of the wedge-shaped region, as can be seen from its formulation and accompanying figures.

First, we place an \emph{origin tile}, a congruent copy of the prototile, starting at the origin~$O$ and in the same orientation.
We refer to the starting point and the base edge of a wedge as those of its origin tile.

Next, we place two additional tiles: one shares its base edge with the opposite edge of the origin tile's base edge, and the other shares its neighbor edge with the opposite edge of the origin tile's neighbor edge.
These three tiles meet at the endpoint of the origin tile and together form the core triangular structure from which the wedge is extended.

\begin{figure}[htbp]
  \centering
  \includegraphics[width=0.45\textwidth]{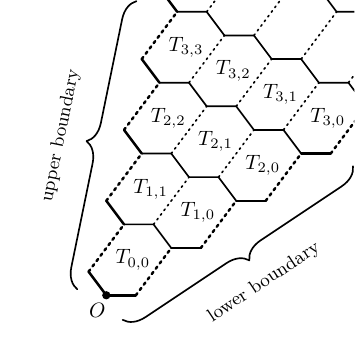}
  \caption{A conceptual diagram of an initial wedge. Each tile is placed such that its base edge or neighbor edge aligns with the opposite edge of a previously placed tile.}
  \label{fig:wedge_concept}
\end{figure}

Formally, we continue placing tiles iteratively along two directions, using the following translation vectors:
\[
\vec{d}_0 = \sum_{j=0}^{k-1} \vec{v}_{\ell_j}, \quad \vec{d}_1 = \vec{v}_k - \vec{v}_0.
\]
Here, $\vec{d}_0$ corresponds to the row shift, and $\vec{d}_1$ to the in-row shift.
For integers $r$ and $c$ such that $0 \le c \le r$, we place a tile $T_{r\!, c}$ starting at $r \vec{d}_0 + c \vec{d}_1$.
This construction is illustrated in Figure~\ref{fig:wedge_concept}.

The resulting wedge is bounded by two polygonal paths extending from the origin: the \emph{lower boundary} and the \emph{upper boundary}.
By construction, the direction sequence of the lower boundary is a periodic repetition of the first~$k$ terms of $(\ell_j)$; similarly, that of the upper boundary follows the same periodic structure defined by $(u_j)$.
This property becomes essential in Section~\ref{sec:proof}, where we formally prove the validity of the construction.

\section{Tiling the Plane}
\label{sec:tiling}

In this section, we describe how to tile the plane using the wedge constructed in the previous section.

We begin with \emph{wedge~$0$}, the initial wedge described in the previous section.
As shown earlier, this wedge covers a region bounded by two infinite polygonal paths extending from the origin.
We denote the direction sequences along the lower and upper boundaries of wedge~$0$ by $(\ell^0_j)$ and $(u^0_j)$, respectively.

We proceed counterclockwise, placing wedges labeled $1, 2, \dots$ in increasing order of direction indices, until a certain condition is satisfied.
To describe the placement rule, we define two boundary paths:
\begin{itemize}
  \item The \emph{base} is the lower boundary of the region, defined by the direction sequence $(b_j) \coloneqq (\ell^0_j)$, which remains fixed throughout the process.
  \item The \emph{front} is the upper boundary of the region covered so far. We denote the front direction sequence after placing wedge~$i$ by $(f^i_j)$. Initially, we have $(f^0_j) = (u^0_j)$, and the front evolves as additional wedges are placed.
\end{itemize}

Each new wedge~$i$ is placed so that its base edge aligns with the first edge in the front that points in direction~$i$.
In other words, we rotate a copy of the initial wedge by $\frac{2\pi i}{n}$ about the origin.
We then identify the minimal index $j^*$ such that $f^{i-1}_{j^*} = i$, and follow the front along the directions $f^{i-1}_0, f^{i-1}_1, \dots, f^{i-1}_{j^*-1}$ to determine the starting point of wedge~$i$, and attach the wedge at that position.
Wedge~$i$ is thus a copy of wedge~$0$ that has been rotated counterclockwise by $\frac{2\pi i}{n}$ about the origin, and then translated by the following vector:
\[
\sum_{j=0}^{j^*-1} \vec{v}_{f^{i-1}_j}.
\]

Accordingly, we denote the direction sequences of the lower and upper boundaries of wedge~$i$ by $(\ell^i_j)$ and $(u^i_j)$, respectively.
Clearly, each of these sequences is obtained by adding $i$ to the corresponding sequence of wedge~$0$:
\[
  \ell^i_j = \ell^0_j + i, \quad
  u^i_j = u^0_j + i.
\]

As we have seen, each wedge is placed by simply aligning its base edge with the first front edge that points in direction~$i$.
Remarkably, its entire lower boundary aligns precisely with the current front, and its upper boundary extends the structure in a consistent and predictable manner---a property we will formally establish in the next section.

We continue placing wedges until the wedge index reaches
\[
w \coloneqq \begin{cases}
k & \text{(without offset)}, \\
\frac{n}{2} = tk - m & \text{(with offset)}.
\end{cases}
\]

After placing wedges counterclockwise as described above, the remainder of the tiling depends on whether the construction involves offset.

\begin{itemize}
  \item \textbf{Without offset} ($n = tk$):
  Placing $k$ wedges yields a region whose front direction sequence $(f^{k-1}_j)$ satisfies
  \[
    f^{k-1}_j = b_j + k.
  \]
  That is, the front coincides with the base rotated by $\frac{2\pi k}{n} = \frac{2\pi}{t}$ about the origin.
  Therefore, placing $t$ copies of this region, each rotated successively by $\frac{2\pi}{t}$, results in a complete tiling of the plane.

  \item \textbf{With offset} ($n = 2(tk - m)$):
  Placing $\frac{n}{2} = tk - m$ wedges yields a region whose front direction sequence $(f^{\frac{n}{2}-1}_j)$ satisfies
  \[
    f^{\frac{n}{2}-1}_j = b_{j+1} + \frac{n}{2}.
  \]
  That is, the front coincides with the truncated base---excluding its first edge---rotated by $\pi$ about the midpoint of the first base edge.
  Therefore, placing a copy of this region, rotated by $180^\circ$ about that midpoint, results in a complete tiling of the plane.
\end{itemize}

\section{Proof of Completeness}
\label{sec:proof}

\subsection{Key Definitions and Observations}

We begin by introducing a fundamental sequence, which we call the \emph{shifted modular progression}.

\begin{definition}[Shifted Modular Progression]
\label{def:shifted_modular_progression}
Let $m$ and $k$ be fixed positive integers.
For each integer $i \ge 0$, define a sequence $(s^i_j)$ as follows:
\[
s^0_j = jm \bmod k, \quad
s^{i+1}_j = \begin{cases}
s^i_j + k & \text{if } s^i_j = i, \\
s^i_j & \text{otherwise}.
\end{cases}
\]
\end{definition}
\noindent
That is, $(s^0_j)$ is a standard modular progression of multiples of $m$ modulo~$k$, and $(s^{i+1}_j)$ is obtained from $(s^i_j)$ by replacing each occurrence of $i$ with $i + k$.

Next, we define a condition that characterizes the behavior of these sequences.

\begin{definition}[Characterization Condition $\textup{C}(i)$]
\label{def:characterization_condition}
A sequence $(a_j)$ is said to satisfy $\textup{C}(i)$ if the following three conditions hold:
\begin{enumerate}
  \item For all integer $j \ge 0$, \ $i \le a_j < i + k$.
  \item The first~$k$ terms $a_0, a_1, \dots, a_{k-1}$ form a permutation of $i, i + 1, \dots, i + k - 1$. In addition, the sequence is periodic with period~$k$.
  \item For all integer $j \ge 0$, the recurrence relation holds:
  \[
    a_{j+1} = \begin{cases}
      a_j + m & \text{if } a_j + m < i + k, \\
      a_j + m - k & \text{if } a_j + m \ge i + k.
    \end{cases}
  \]
\end{enumerate}
\end{definition}

From Definition~\ref{def:characterization_condition}, we get the following:

\begin{observation}
\label{obs:preservation-of-c}
If a sequence $(a_j)$ satisfies $\textup{C}(i)$, then the sequence obtained by removing its first term also satisfies $\textup{C}(i)$.
Furthermore, adding a fixed integer $c$ to each term of $(a_j)$ results in a sequence that satisfies $\textup{C}(i + c)$.
\end{observation}

The following holds by induction on $i$, applied to Definition~\ref{def:shifted_modular_progression}:

\begin{observation}
\label{obs:c-for-s}
The sequence $(s^i_j)$ satisfies $\textup{C}(i)$ for all integer $i \ge 0$.
\end{observation}

We also note the following, which can be verified by induction on $i$: each replacement step adds $k$ to some entries of the sequence.

\begin{observation}
\label{obs:s0-and-si}
For all integers $i, j \ge 0$, there exists an integer $a \ge 0$ such that $s^i_j = s^0_j + ak$.
\end{observation}

Finally, we recall the construction of $(b_j) = (\ell^0_j)$, and observe that replacing each occurrence of $i$ in $(\ell^i_j)$ with $i + k$ gives the sequence $(u^i_j)$.

\begin{observation}
\label{obs:b-ell-and-u}
For all integer $j \ge 0$, \ $b_j = \ell^0_j = jm \bmod k = s^0_j$.
Also, $(u^i_j)$ is obtained from $(\ell^i_j)$ by replacing each occurrence of $i$ with $i + k$.
\end{observation}

We now verify that the sequence $(\ell^i_j)$ also satisfies $\textup{C}(i)$.

\begin{lemma}
\label{lem:c-for-ell}
The sequence $(\ell^i_j)$ satisfies $\textup{C}(i)$ for all integer $i \ge 0$.
\end{lemma}

\begin{proof}
By Observation~\ref{obs:b-ell-and-u}, $(\ell^0_j) = (s^0_j)$ and thus satisfies $\textup{C}(0)$.
For all integer $i \ge 1$, we have $\ell^i_j = \ell^0_j + i$.
Since $(\ell^0_j)$ satisfies $\textup{C}(0)$, it follows from Observation~\ref{obs:preservation-of-c} that $(\ell^i_j)$ satisfies $\textup{C}(i)$.
\end{proof}

\subsection{Main Proof}

\begin{proposition}[Validity of the Procedure]
\label{prop:validity-procedure}
For each $0 \le i < w$, the procedure described in Section~\ref{sec:tiling} is well-defined for placing wedge~$i$.
Moreover, if $i \ge 1$, the lower boundary of wedge~$i$ aligns exactly with the front after placing wedge~$i-1$.
Finally, the front after placing wedge~$i$ is given by $(f^i_j) = (s^{i+1}_j)$.
\end{proposition}

\begin{proof}
We proceed by induction on $i$.

For the base case $i = 0$, wedge~$0$ is constructed with its starting point at the origin.
The front direction sequence is given by $(f^0_j) = (u^0_j)$.
By Observation~\ref{obs:b-ell-and-u}, $(u^0_j)$ is obtained from $(\ell^0_j)$ by replacing each occurrence of $0$ with $k$.
Since $(\ell^0_j) = (s^0_j)$, it follows that $(f^0_j) = (s^1_j)$ by Definition~\ref{def:shifted_modular_progression}.

Now suppose that the claim holds for $i-1$, i.e., $(f^{i-1}_j) = (s^i_j)$, and that wedge~$i$ is placed according to the rule described in Section~\ref{sec:tiling}.

By Observation~\ref{obs:c-for-s}, the sequence $(s^i_j)$ satisfies $\textup{C}(i)$.
In particular, its first $k$ terms form a permutation of the integers from $i$ to $i + k - 1$, so there always exists a minimal index $j^*$ such that $f^{i-1}_{j^*} = i$ holds\footnote{Since we have $f^{i-1}_j \equiv s^i_j \equiv jm \pmod{k}$, solving $f^{i-1}_j \equiv i \pmod{k}$ gives $j \equiv im^{-1} \pmod{k}$. Thus, we can define $j^* \coloneqq im^{-1} \bmod k$, where $m^{-1}$ denotes the modular inverse of $m$ modulo~$k$.}.
The direction sequence of the lower boundary of wedge~$i$ is $(\ell^i_j)$, which also satisfies $\textup{C}(i)$ by Lemma~\ref{lem:c-for-ell}.
Since both sequences satisfy the same recurrence relation (Condition~3 in Definition~\ref{def:characterization_condition}), their values coincide from the matching index onward.
Therefore, the lower boundary of wedge~$i$ aligns exactly with the previous front.

Next, we analyze how the front direction sequence changes upon placing wedge~$i$.
The portion of the front direction sequence that matches $(\ell^i_j)$ is replaced with $(u^i_j)$.
By Observation~\ref{obs:b-ell-and-u}, the sequence $(u^i_j)$ is obtained from $(\ell^i_j)$ by replacing each occurrence of $i$ with $i + k$.
Since we align the base edge of wedge~$i$ to the \emph{first} edge in the current front that points in direction~$i$, every occurrence of $i$ in $(s^i_j)$ is replaced with $i + k$, yielding $(s^{i+1}_j)$ by Definition~\ref{def:shifted_modular_progression}.

This completes the induction.
\end{proof}

\begin{proposition}[Symmetry in the Case Without Offset]
\label{prop:symmetry-no-offset}
After placing $w = k = \frac{n}{t}$ wedges as described in Section~\ref{sec:tiling}, the resulting region has a front direction sequence $(f^{k-1}_j)$ satisfying
\[
f^{k-1}_j = b_j + k.
\]
That is, the front coincides with the base rotated by $\frac{2\pi k}{n} = \frac{2\pi}{t}$ about the origin.
Therefore, placing $t$ copies of this region, each rotated successively by $\frac{2\pi}{t}$, results in a complete tiling of the plane.
\end{proposition}

\begin{proof}
By Proposition~\ref{prop:validity-procedure}, after placing $w = k = \frac{n}{t}$ wedges, the front direction sequence is $(f^{k-1}_j) = (s^k_j)$, which satisfies $\textup{C}(k)$ by Observation~\ref{obs:c-for-s}.
By Definition~\ref{def:shifted_modular_progression}, it follows that $0 \le s^0_j < k$ for all integer $j \ge 0$.
On the other hand, since $(s^k_j)$ satisfies $\textup{C}(k)$, we have $k \le s^k_j < 2k$ for all integer $j \ge 0$.
Moreover, by Observation~\ref{obs:s0-and-si}, each $s^k_j$ can be written as $s^k_j = s^0_j + a_j\cdot k$ for some integer $a_j \ge 0$.
Given the ranges above, we see that $a_j = 1$ for all integer $j \ge 0$, and thus $s^k_j = s^0_j + k$.
Therefore, we have $f^{k-1}_j = s^k_j = s^0_j + k = b_j + k$ for all integer $j \ge 0$.

Since $\vec{v}_i$ is defined as a unit vector with angle $\frac{2\pi i}{n}$, adding $k = \frac{n}{t}$ to a direction corresponds to a rotation by $\frac{2\pi\cdot\frac{n}{t}}{n} = \frac{2\pi}{t}$.
Thus, the front coincides with the base rotated by $\frac{2\pi}{t}$ about the origin.

This relation guarantees that placing $t$ copies of the region, each rotated successively by $\frac{2\pi}{t}$, completes the tiling without gaps or overlaps.
\end{proof}

\begin{proposition}[Symmetry in the Case With Offset]
\label{prop:symmetry-with-offset}
After placing $w = \frac{n}{2}$ wedges as described in Section~\ref{sec:tiling}, the resulting region has a front direction sequence $(f^{\frac{n}{2}-1}_j)$ satisfying
\[
f^{\frac{n}{2}-1}_j = b_{j+1} + \frac{n}{2}.
\]
That is, the front coincides with the truncated base---excluding its first edge---rotated by $\pi$ about the midpoint of the first base edge.
Therefore, placing a copy of this region, rotated by $180^\circ$ about that midpoint, results in a complete tiling of the plane.
\end{proposition}

\begin{proof}
By Proposition~\ref{prop:validity-procedure}, after placing $w = \frac{n}{2} = tk - m$ wedges, the front direction sequence is $(f^{w-1}_j) = (s^w_j)$, which satisfies $\textup{C}(w)$ by Observation~\ref{obs:c-for-s}.
To determine the starting direction of the front, we examine $s^w_0$, the first term of $(s^w_j)$.
By Observation~\ref{obs:s0-and-si}, $s^w_0$ is a multiple of $k$.
Moreover, since $(s^w_j)$ satisfies $\textup{C}(w)$, $s^w_0$ must satisfy
\[
tk - m \le s^w_0 < (t+1)k - m.
\]
The only multiple of $k$ in this interval is $tk$.
Hence, the first term of the front direction sequence is $tk = \frac{n}{2} + m$, and the full sequence satisfies $\textup{C}(\frac{n}{2})$.

Now consider the base direction sequence $(b_j)$, which is defined as $(\ell^0_j)$.
We remove its first term to obtain the sequence $(b_{j+1})$, which begins with $m$.
Since $(\ell^0_j)$ satisfies $\textup{C}(0)$ by Lemma~\ref{lem:c-for-ell}, the truncated sequence $(b_{j+1})$ also satisfies $\textup{C}(0)$ by Observation~\ref{obs:preservation-of-c}.
By adding $\frac{n}{2}$ to each term of $(b_{j+1})$, we obtain a sequence with first term $\frac{n}{2} + m$ satisfying $\textup{C}(\frac{n}{2})$ (again by Observation~\ref{obs:preservation-of-c}).
Therefore, the front direction sequence coincides with the truncated and shifted base direction sequence.

Since $\vec{v}_i$ is defined as a unit vector with angle $\frac{2\pi i}{n}$, adding $\frac{n}{2}$ to a direction corresponds to a rotation by $\frac{2\pi\cdot\frac{n}{2}}{n} = \pi$.
Thus, the front coincides with the truncated base rotated by $\pi$ about the midpoint of the first base edge.

This relation guarantees that placing a copy of the region, rotated by $180^\circ$ about that midpoint, completes the tiling without gaps or overlaps.
\end{proof}

\begin{theorem}[Completeness of Modulo Krinkle Tiling]
The tiling constructed as described in Section~\ref{sec:tiling} correctly fills the plane.
\end{theorem}

\begin{proof}
Immediate from Propositions~\ref{prop:symmetry-no-offset} and~\ref{prop:symmetry-with-offset}, which together cover all cases.
\end{proof}

In what follows, we refer to the tiling constructed in Section~\ref{sec:tiling} simply as the \emph{$(m, k, n)$-Modulo Krinkle tiling}.
This notation is justified because, in most cases, the parameters $t$ and whether the construction involves an offset can be uniquely determined from the triple $(m, k, n)$.
For instance, when $n$ is a multiple of $k$, we are in the case without offset, and $t = \frac{n}{k}$.
Otherwise, we are in the offset case, and $t$ can be recovered from the relation $n = 2(tk - m)$.

\clearpage
\section{Modulo-Staggered Rotational Symmetry}
\label{sec:symmetry}

The construction of a wedge clearly demonstrates that our prototile admits a periodic tiling of the plane, and thus is not an aperiodic monotile.
After all, if the sole objective is to produce a non-periodic tiling, one can easily achieve this by modifying a regular tiling using elementary shapes such as rectangles or $45$-$45$-$90$ right triangles.
Does this imply that our tiling system offers nothing more than these trivial examples?

We argue otherwise.
In fact, the resulting tilings exhibit distinct patterns, suggesting a new kind of structural regularity not captured by conventional frameworks.

From the perspective of Euclidean symmetry, a Modulo Krinkle tiling without offset exhibits $t$-fold rotational symmetry, and one with offset exhibits only point symmetry.
Nevertheless, we observe that each of these tilings is composed of $n$ congruent wedge-shaped regions, rotated successively by $\frac{2\pi}{n}$ around a center.
This raises the question: can we regard this structure as exhibiting a kind of pseudo-$n$-fold rotational symmetry?

We now relax or extend the concept of symmetry: we may consider figures that are divided into multiple congruent parts, such that any permutation of the parts preserves the overall structure.
Such a perspective---based on repetition of structure across congruent regions---enables a broader interpretation of symmetry.
In addition, the fact that our tilings place $n$ congruent regions at equal angular steps makes them resemble standard $n$-fold rotational symmetry even more closely\footnote{Notably, dividing the plane into $n$ congruent wedge-shaped regions using curves defined by shifted modular progressions itself appears to be a distinctive construction, not commonly found in the literature.}.

In this spirit, we propose a generalized notion of rotational symmetry.
We~say that a figure exhibits \emph{$n$-fold $(m, k)$-modulo-staggered rotational symmetry} if it is divided into $n$ wedge-shaped regions according to a $(m, k, n)$-Modulo Krinkle tiling, with each region exhibiting the same internal structure.

This concept offers a way to formally capture the distinct pattern observed in Modulo Krinkle tilings, positioning them as meaningful instances of this extended notion of symmetry.

\clearpage
\section{Conclusion}

We have introduced a family of non-periodic tilings that are defined and verified using only elementary tools, yet exhibit a previously undocumented form of structural regularity. Future work may include further characterization of these tilings and deeper interpretations of their structure.

Beyond mathematical interest, the construction lends itself to a variety of applications.
Its distinct aesthetic properties make it a natural candidate for generative design.
In addition, its elementary construction makes it suitable for educational purposes.
Moreover, some variants are structurally simple enough to be implemented in real-world materials, such as architectural tiling.
In particular, the capability of placing tiles with gradually changing orientation may prove useful in practical settings such as road paving or other surface applications.
These aspects, though not explored in this paper, suggest a broader potential for the proposed system.

\bibliographystyle{unsrt}
\bibliography{references}

\appendix

\section{Gallery of Generated Patterns}

The following figures illustrate diverse instances of the Modulo Krinkle tiling.
These patterns are generated with varying parameters and highlight the visual diversity and design potential of the system.

\begin{figure}[htbp]
  \centering

  \begin{subfigure}[b]{0.45\textwidth}
    \centering
    \includegraphics[width=\textwidth]{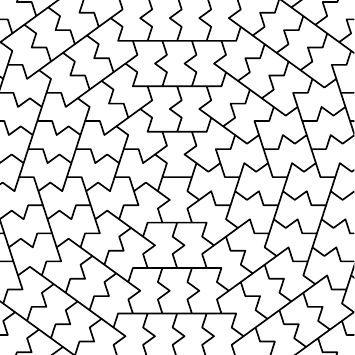}
    \caption{$(m, k, n) = (2, 5, 10)$}
    \label{fig:gallery_a}
  \end{subfigure}
  \hspace{0.03\textwidth}
  \begin{subfigure}[b]{0.45\textwidth}
    \centering
    \includegraphics[width=\textwidth]{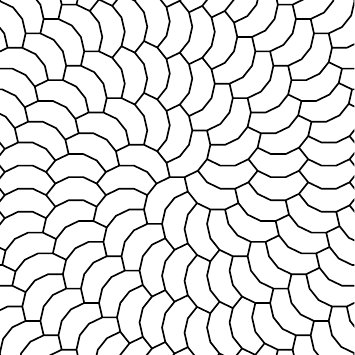}
    \caption{$(m, k, n) = (1, 4, 14)$}
    \label{fig:gallery_b}
  \end{subfigure}

  \vspace{0.02\textwidth}

  \begin{subfigure}[b]{0.45\textwidth}
    \centering
    \includegraphics[width=\textwidth]{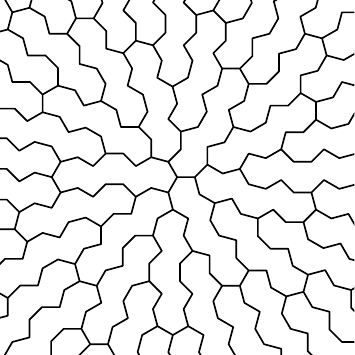}
    \caption{$(m, k, n) = (3, 8, 24)$}
    \label{fig:gallery_c}
  \end{subfigure}
  \hspace{0.03\textwidth}
  \begin{subfigure}[b]{0.45\textwidth}
    \centering
    \includegraphics[width=\textwidth]{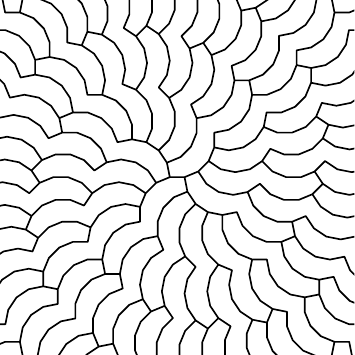}
    \caption{$(m, k, n) = (2, 9, 32)$}
    \label{fig:gallery_d}
  \end{subfigure}

  \vspace{0.02\textwidth}

  \begin{subfigure}[b]{0.45\textwidth}
    \centering
    \includegraphics[width=\textwidth]{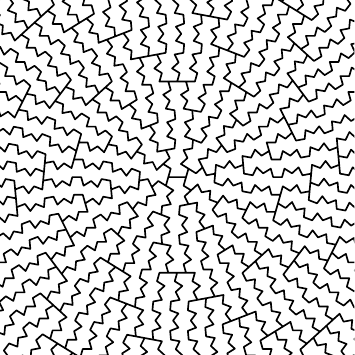}
    \caption{$(m, k, n) = (7, 17, 34)$}
    \label{fig:gallery_e}
  \end{subfigure}
  \hspace{0.03\textwidth}
  \begin{subfigure}[b]{0.45\textwidth}
    \centering
    \includegraphics[width=\textwidth]{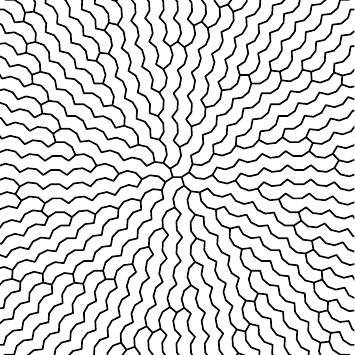}
    \caption{$(m, k, n) = (5, 16, 54)$}
    \label{fig:gallery_f}
  \end{subfigure}

  \caption{A variety of $(m, k, n)$-Modulo Krinkle tilings.}
  \label{fig:gallery}
\end{figure}

The pattern shown in Figure~\ref{fig:gallery_b} was originally presented by Jan Sallmann-Räder in a social media post, where regular heptagons were manually assembled.
Our construction significantly generalizes this example, providing an extensible system that produces infinitely many non-periodic patterns.

\end{document}